\documentclass[reqno]{amsart}
\usepackage{amssymb}  
\usepackage{amsmath} 
\usepackage{amsthm}  
\usepackage{stmaryrd} 

\newtheorem{thm}{Theorem}
\newtheorem{prop}[thm]{Proposition}

\renewcommand{\Im}{\operatorname{Im}}
\newcommand{\Ker}{\operatorname{Ker}}
\newcommand{\nor}{\trianglelefteqslant}

\title{On the embedding of central extensions into wreath products}
\author{Andrei V. Zavarnitsine}

\address{\textup{\scriptsize
Andrei V. Zavarnitsine\\
Sobolev Institute of Mathematics\\
4, Koptyug av.\\
630090, Novosibirsk, Russia\\
}}
\email{zav@math.nsc.ru}

\date{}

\begin{document}
\begin{abstract} We find a necessary condition for the embedding of a central extension of a group $G$ with
elementary abelian kernel into the wreath product that corresponds to a permutation action of $G$. The proof uses purely group-theoretic methods.



\end{abstract}
\maketitle

The finite group $G=\operatorname{PSL}_2(q)$, $q$ odd, acts naturally by permutations on the projective line of order $q+1$. In
\cite{13Zav}, we studied the embedding of $\operatorname{SL}_2(q)$ into the wreath product of $\mathbb{Z}/2\mathbb{Z}$ and $G$ that
corresponds to this permutation action. This problem was generalized in \cite{16Zav} to arbitrary groups
$\operatorname{PSL}_n(q)$ and their central extensions with kernel of prime order. In this paper, we obtain a further
generalisation of these results.

Let $G$ be a finite group, $\Omega$ a finite set, and let $\rho: G\to \operatorname{Sym}(\Omega)$ be a permutation
representation. For $\omega\in \Omega$ we denote by $\operatorname{St}(\omega)$ the stabilizer of $\omega$ in $G$.

For a commutative unital ring $A$ of prime characteristic $p$, consider the (right) permutation $AG$-module $V$ corresponding to $\rho$
with basis (identified with) $\Omega$ and its submodule
\begin{equation}\label{iv}
0\to I\to V
\end{equation}
generated by $\omega_0=\sum_{\omega\in \Omega}\omega$. Clearly, there is an $AG$-module isomorphism $\alpha:A\to I$ defined by
$\alpha:1\mapsto \omega_0$. Let $G\rightthreetimes V$ denote the natural semidirect product.

Assume that we also have a central extension
\begin{equation}\label{ce}
1\to A \stackrel{\iota}{\to} H\stackrel{\pi}{\to} G\to 1,
\end{equation}
i.\,e. one with $\Im \iota\leqslant \mathrm{Z}(H)$, where we identify $A$ with its additive group $A^+$. A subgroup $S\leqslant
G$ is {\em liftable} to $H$ if the full preimage $S\pi^{-1}$ splits over $\Im \iota$. We say that $H$ is a {\em subextension} of
$G\rightthreetimes V$ with respect to the embedding (\ref{iv}), if there is an embedding $\beta: H\to G\rightthreetimes V$ such
that the following diagram commutes
\begin{equation}\label{sube}
\begin{array}{cc@{}c@{}c@{\ }c@{\,}cccl}
1&\longrightarrow&A         &\stackrel{\iota}{\longrightarrow}&H         &\stackrel{\pi}{\longrightarrow}&G&\longrightarrow&1\\[2pt]
 &   &\phantom{\alpha}\downarrow\alpha&   &\phantom{\beta}\downarrow\beta&   &\|&&   \\[2pt]
1&\longrightarrow&V         &\longrightarrow&G\rightthreetimes V         &\longrightarrow&G&\longrightarrow&1,
\end{array}
\end{equation}
where we identify $I=\operatorname{Im}\alpha$ with its image in $V$ under (\ref{iv}).

The main result to be proved in Section \ref{nc} is the following necessary condition.

\begin{thm}\label{main} In the above notation, if  a central extension $H$ is a subextension
of $G\rightthreetimes V$ with respect to the embedding $(\ref{iv})$ then $\operatorname{St}(\omega)$ is liftable to $H$ for every
$\omega\in \Omega$.
\end{thm}

The proof generalises some ideas presented in \cite{13Zav,16Zav}. In particular, we also prove an auxiliary result about
presentations of $p$-groups.

Let $F=F\langle X\rangle$ be a free group with basis $X$. Every $w\in F$ can be written in the form
$$
w=x_1^{\varepsilon_1}\ldots x_t^{\varepsilon_t},
$$
where $x_i\in X$ and $\varepsilon=\pm 1$. For $x\in X$, we define
$$
\mu_x(w)=\sum_{x_i=x} \varepsilon_i.
$$

The following fact is proven in Section \ref{pprg}.

\begin{prop}\label{pgr} Every finite $p$-group $P$ has a finite presentation $\langle\,X\mid R\,\rangle$ such that $\mu_x(r)\equiv 0\pmod p$ for all
$x\in X$ and $r\in R$.
\end{prop}

\section{Fox derivatives}

Let $X=\{x_1,\ldots,x_n\}$ and let $F=F\langle X\rangle$ be a free group with basis $X$. Recall that the (right) {\em Fox
derivative} $\partial/\partial x_i:F\to\mathbb{Z}F$ is the map satisfying $\partial x_j/\partial x_i=\delta_{ij}$, $1\leqslant j
\leqslant n$ and
$$
\frac{\partial(uv)}{\partial x_i}= \frac{\partial u}{\partial x_i}\, v + \frac{\partial v}{\partial x_i}
$$
for all $u,v\in F$ and $1\leqslant i \leqslant n$. Let $w=w(x_1,\ldots,x_n)\in F$ and write
$$
w=x_{i_1}^{\varepsilon_1}\ldots x_{i_l}^{\varepsilon_l},
$$
where $x_{i_k}\in X$ and $\varepsilon_k=\pm 1$ for all $k$. It can be shown \cite[Proposition 2.73]{05HolEicObr} that
$$
\frac{\partial w}{\partial x_i}=\sum_{\{ k\, \mid\, i_k=i \}} \varepsilon_k f_k,
$$
where
\begin{equation}\label{fk}
f_k=\left\{
\begin{array}{rl}
  x_{i_{k+1}}^{\varepsilon_{k+1}}\ldots x_{i_l}^{\varepsilon_l}, & \varepsilon_k=1, \\
  x_{i_k}^{\varepsilon_k}x_{i_{k+1}}^{\varepsilon_{k+1}}\ldots x_{i_l}^{\varepsilon_l}, & \varepsilon_k=-1.
\end{array}
\right.
\end{equation}

Let $G$ be a group and $V$ a $G$-module. Fixing a homomorphism
$F\to G\rightthreetimes V$ we write the image of each $x_i$ as $g_iv_i$ for suitable $g_i\in G$, $v_i\in V$. Then using the
additive notation in $V$ we can write
\begin{equation}\label{wg1n}
w(g_1v_1,\ldots,g_nv_n)=w(g_1,\ldots,g_n)(v_1\frac{\partial w}{\partial g_1}+\ldots+v_n\frac{\partial w}{\partial g_n}),
\end{equation}
where $\partial/\partial g_i$ is the short-hand notation for the composition of  $\partial/\partial x_i$ and the homomorphism
$\mathbb{Z}F\to\mathbb{Z}G$ which extends the map $x_i\mapsto g_i$, $i=1,\ldots,n$. For details, see  \cite[\S 1.9]{06Kuz}.

\section{Presentations of group extensions}

Let
$$
1\to N\stackrel{\iota}{\to}G\stackrel{\pi}{\to} Q\to 1
$$
be a short exact sequence of groups.  Suppose that $N$ has presentation $\langle\, \overline Y\mid \overline S \,\rangle$ and
$Q$ has presentation $\langle\, \overline X \mid \overline R \,\rangle$. Using this information it is possible to describe a presentation of $G$.
Let $Y$ be the image of $\overline Y$ under $\iota: \overline y\mapsto y$ and let
$$
S=\{\,s\mid \overline s \in \overline S\,\},
$$
where $s$
is the word in $Y$ obtained from $\overline s$ by replacing each $\overline y$ with $y$. Choose $X\subseteq G$ so that
$x\pi = \overline x \in \overline X $ for all $x\in X$. For every $\overline r \in \overline R$, let $r$
be the word in $X$ obtained from $\overline r$ by replacing each $\overline x$ with $x$. Clearly, $r$
as an element of $G$ lies in $\Ker\pi = \Im \iota$ and so is a word, say $w_r$, in $Y$. Define
$$
R=\{\,rw_r^{-1}\mid \overline r \in \overline R\,\}.
$$
Also, since $\Im \iota \nor G$, the element $x^{-1}yx$ lies in $\Im \iota$ for all $y\in Y, x\in X$ and so is a word, say $w_{xy}$, in $Y$.
We set
$$
T=\{\,x^{-1}yxw_{xy}^{-1}\mid \overline x \in \overline X,\  \overline y \in \overline Y\ \}.
$$

\begin{prop}\cite[Proposition 10.2.1]{97Joh}, \cite[Proposition 2.55]{05HolEicObr} \label{ext_pres} In the above notation,
\begin{equation}\label{pres}
\langle\, X \cup Y\mid R \cup S \cup T \,\rangle
\end{equation}
is a presentation of $G$.
\end{prop}

\section{Proof of Proposition \ref{pgr} \label{pprg}}

Recall that $\Omega_1(P)$ denotes the subgroup of a $p$-group $P$ generated by all elements of order~$p$.

\begin{proof} We use induction on $|P|$. If $|P|=1$, the claim holds. Assume $|P|>1$ and let $N=\Omega_1(\operatorname{Z}(P))$. Note that $N$ is a nontrivial
elementary abelian $p$-group and
$$
1 \to N \stackrel{\iota}{\to} P \to Q\to 1
$$
is a central extension. By induction, $Q$ has a finite presentation $\langle\, \overline X \mid \overline R \,\rangle$ that
satisfies the required property. Clearly, $N$ also has a  presentation $\langle \overline Y \mid \overline S \,\rangle$, where
$\overline Y$ is finite and
$$
\overline S = \{\,\overline y^{\,p}, \, [\overline y_1,\overline y_2]\ \mid \ \overline y ,\overline y_1,\overline y_2\in \overline Y \,\},
$$
which has the required property. Note that we may take any basis of $N$ as $\overline Y$. We define the sets of generators $X$
and $Y$ and relators $R$, $S$, and $T$ as before Proposition~\ref{ext_pres}, where $G=P$. Since the relators in $S$ are rewritten
from those of $\overline S$, they have the required property, i.\,e., the exponent sum for each generator in each relator is a
multiple of $p$. Also, since $\Im \iota$ is central in $P$, we have $w_{xy}=y$ for all $x\in X$, $y\in Y$, and so $T$ consists of
commutators which have the required property. We now consider the relators $rw_r^{-1}$ in $R$. Some of them will be eliminated,
while in others we will replace the subwords $w_r$ with ones satisfying the required property.

Indeed, we can choose a maximal linearly independent subset of
$$
W=\{\,w_r\mid r\in R\,\}\subseteq\Im\iota
$$
and complete it to a basis of $\Im \iota$. As we have mentioned, without loss of generality we may assume that this basis
coincides with $Y$. All  generators $y=w_r\in W\cap Y$ may be eliminated, because we have a relation $w_r=r$ and $r$ does not
involve any $y\in Y$. The remaining words $w_r\in W\setminus Y$ are linear combinations of such generators, hence after the elimination
they will become words in $R$ which satisfy the needed property by induction. The words in $S\cup T$ are commutators and powers $y^p$, hence will
retain the needed property, too. The resulting presentation of $P$ clearly has the required property.
\end{proof}

\section{Proof of main theorem \label{nc}}

The following result will be used.

\begin{prop}[Gasch\"utz' Theorem {\cite[(10.4)]{00Asc}}]\label{gth}
Let $p$ be a prime, $V$ a normal abelian $p$-subgroup of a finite group $G$, and $P\in \operatorname{Syl}_p(G)$. Then $G$ splits
over $V$ if and only if $P$ splits over $V$.
\end{prop}

We are now ready to prove Theorem \ref{main}.

\begin{proof} We denote $\omega_0=\sum_{\omega\in \Omega}\omega$.
Assume to the contrary that there is $\omega\in \Omega$ such that $S=\operatorname{St}(\omega)$ is not liftable to $H$. Let $P\in
\operatorname{Syl}_p(S)$. Since $A$ is an abelian $p$-group, Proposition \ref{gth} implies that $P$ is not liftable to $H$. Let
$\langle\,X\mid R\,\rangle$ be a finite presentation for $P$ with the property that $\mu_x(r)\equiv 0\pmod p$ for every $x\in X$, $r\in
R$. Such a presentation exists by Proposition \ref{pgr}.

Let $F=F\langle X\rangle$ be the free group with basis $X=\{x_1,\ldots
x_n\}$. For every $x\in X$, we denote $\underline x=x\gamma\in P$,
where $\gamma: F\to P$ is the presentation homomorphism, and  choose $\overline x\in H$ so that
\begin{equation}\label{oxux}
  \overline x\pi = \underline x.
\end{equation}
There exists a relator $r=r(x_1,\ldots
x_n)\in R$ such that $\overline{r}=r(\overline{x}_1,\ldots,\overline{x}_n)\ne 1$ in $H$. Indeed, otherwise the subgroup
$$
\overline{P}=\langle\,\overline x \mid x\in X\,\rangle \leqslant H
$$
would satisfy the same relations as $P$ and so the map $[x \mapsto \overline x]$, $x\in X$, would give rise to a homomorphism
$\sigma: P\to\overline P$ with the property $\sigma\pi=\operatorname{id}_P$. But this means that $\overline{P}$ would be a
splitting of the full preimage $P\pi^{-1}$ contrary to the assumption.

Since
\begin{equation}\label{rux}
\overline{r}\pi=r(\underline{x}_1,\ldots,\underline{x}_n)=1,
\end{equation}
we see that $\overline{r}=a\iota$ for a nonzero $a\in A$. By assumption, $H$ is a subextension of $G\rightthreetimes V$ with
respect to (\ref{iv}). Hence $\overline{r}\beta=a\omega_0$, where the embedding $\beta$ is as in (\ref{sube}). Also, we can write
$\overline x_i\beta =g_i v_i$, $i=1,\ldots,n$, for suitable $g_i\in G$, $v_i\in V$. Observe that $g_i=\underline{x}_i$ due to
(\ref{oxux}) and the commutativity of diagram (\ref{sube}). Let $r=x_{i_1}^{\varepsilon_1}\ldots x_{i_l}^{\varepsilon_l}$ with
$i_k\in \{1,\ldots,n\}$, $\varepsilon_k=\pm 1$, $k=1,\ldots,l$. Define a homomorphism $F\to G\rightthreetimes V$ by extending the
map $x_i\mapsto \underline{x}_i v_i$, $i=1,\ldots,n$. Using (\ref{wg1n}) and (\ref{rux}), we have

\begin{multline}\label{eqa}
 a\omega_0 =\overline{r}\beta = r(\overline{x}_1\beta,\ldots,\overline{x}_n\beta) = r(\underline{x}_1v_1,\ldots,\underline{x}_nv_n) =
r(\underline{x}_1,\ldots,\underline{x}_n)    \\
\times( v_1\frac{\partial r}{\partial \underline{x}_1}+\ldots+v_n\frac{\partial r}{\partial \underline{x}_n} )=v_1\sum_{\{k\,\mid\, i_k=1\}} \varepsilon_k \underline{f}_{\,k}+\ldots + v_n\sum_{\{k\,\mid\, i_k=n\}} \varepsilon_k \underline{f}_{\,k},
\end{multline}
where $f_k\in F$ is given by (\ref{fk}) and $\underline{f}_{\,k}=f_k\gamma \in P$. We can decompose
$$V=A\omega\oplus V_0,$$
where $V_0$ the $A$-linear span of $\Omega\setminus\omega$, and write $v_i=a_i\omega+w_i$, $i=1,\ldots,n$, for suitable $a_i\in
A$ and $w_i\in V_0$. Since $\underline{f}_{\,k}\in S$ stabilizes $\omega$, it also stabilizes $V_0$. Therefore, the right-hand
side of (\ref{eqa}) can be rewritten as
\begin{equation}\label{eqb}
a_1\left(\sum_{\{k\,\mid\, i_k=1\}} \varepsilon_k\right)\omega + w_1' +\ldots + a_n\left(\sum_{\{k\,\mid\, i_k=n\}} \varepsilon_k\right)\omega + w_n',
\end{equation}
where $w_i'=\sum_k \varepsilon_k w_i \underline{f}_{\,k}$ lies in $V_0$ for each $i$. Observe that
$$
\sum_{\{k\,\mid\, i_k=i\}} \varepsilon_k = \mu_{x_i}(r)\equiv 0 \pmod p
$$
for every $i$ by assumption. Since $A$ has characteristic $p$, (\ref{eqb}) equals $\sum_i w_i'=w'$, an element of $V_0$. We now
compare the coefficients of $\omega$ for $w'$ and $a\omega_0$. Since $V$ is free as an $A$-module, these coefficients must
coincide. However, the former is $0$ and the latter is $a\ne 0$, a contradiction.
\end{proof}

\end{document}